\documentclass[twoside,a4paper]{article}
\usepackage{graphicx}
\usepackage[sumlimits,intlimits,reqno]{amsmath}
\usepackage[english]{babel}
\usepackage{amssymb}
\usepackage{amsthm}
\newtheorem*{theorem}{Theorem}
\newtheorem{definition}{Definition}[section]
\newtheorem{lemma}{Lemma}[section]
\numberwithin{equation}{section}

\newcommand\pop[1]{{\partial\over\partial{#1}}}
\newcommand\DoP[2]{{D{#2}\over\partial{#1}}}

\newcommand\II{\operatorname{II}}
\newcommand\diag{\operatorname{diag}}
\newcommand\OB{{}\,\overline{\kern-.4mm B\kern-.3mm}\,}
\newcommand\hA{{}\kern1mm\widehat{\kern-1mm A}}
\newcommand\hB{{}\,\widehat{\kern-.3mm B}}
\newcommand\olr{{}\,\overline{\!R}}

\begin{document}

\title{The Sectional Curvature Remains Positive when Taking Quotients
by Certain Nonfree Actions%
\thanks{The author was supported by
the Russian Foundation for Basic Research (Grant~06-01-0094a).}}
\author{Semyon Dyatlov}
\date{October 2007}
\markboth{Semyon Dyatlov}{The Sectional Curvature Remains Positive when Taking Quotients}

\maketitle

\begin{abstract}
\noindent We study some cases when the sectional curvature remains
positive under the taking of quotients by certain nonfree isometric
actions of Lie groups. We consider the actions of the groups~$S^1$
and~$S^3$ such that the quotient space can be endowed with a smooth
structure using the fibrations~$S^3/S^1{\simeq}S^2$ and
$S^7/S^3{\simeq}\nobreak S^4$. We prove that the quotient space
carries a metric of positive sectional curvature, provided that the
original metric has positive sectional curvature on all 2-planes
orthogonal to the orbits of the action.
\end{abstract}

\section*{Introduction}

In this paper we consider a certain method for constructing closed
Riemannian manifolds of positive sectional curvature. There are only
few known examples of such manifolds (one can find a list of all known
examples, for example, in the introduction to the paper~\cite{Ba}).

The main tool for constructing new examples is the taking of quotients
by free isometric actions of Lie groups. More
precisely, suppose that a Lie group~$G$ acts on a closed smooth
Riemannian manifold~$M$ by isometries; i.e., the metric of~$M$ is
preserved by the action. Suppose further that this action is free;
i.e., for every point $p\in M$ and every element $g\in G$ such
that~$g\neq 1$ we have $g\cdot p\neq p$. In this case the quotient
space~$M/G$ has the structure of a smooth manifold such that the
projection map~$\tau:M\to M/G$ is smooth.

Introduce a metric on~$M/G$ as follows. Let~$X$ and $Y$ be two vectors
tangent to~$M/G$ at the point~$q$, and let $q=\tau(p)$ for some $p\in
M$. Then there exist unique vectors~$X^\prime$ and~$Y^\prime$ tangent
to~$M$ at the point~$p$ that are orthogonal to the orbit of~$G$ and
mapped by~$\tau$ to~$X$ and~$Y$, respectively. Let $g$ be the metric
of~$M$. Put $g_*(X,Y)=g(X^\prime,Y^\prime)$. Since~$G$ acts by
isometries, this relation defines correctly some metric~$g_*$
on~$M/G$, which we call the {\it quotient metric\/} of the metric~$g$.

To compute the curvature of the quotient metric, one uses
the curvature formulas for Riemannian submersions. A map~$f:M\to N$
between manifolds is a~{\it submersion} iff it is surjective together
with its differential at each point. Let $f$ be a submersion between
Riemannian manifolds $(M, g_M)$ and $(N, g_N)$. Consider two
vectors~$X^\prime$ and~$Y^\prime$ at a point~$p\in M$ that are
orthogonal to the submanifold $f^{-1}\big(f(p)\big)$. Suppose that~$f$
maps them to the respective vectors~$X$ and $Y$. The map~$f$ is a {\it
Riemannian submersion} iff $g_M(X^\prime,Y^\prime)=g_N(X,Y)$ for
all~$X^\prime$ and~$Y^\prime$. If~$f$ is a Riemannian submersion and
$R_M$ and $R_N$ are the curvature tensors of the manifolds~$M$
and~$N$, then the relations on the curvature found in~\cite{O'Ne}
imply that
$$
R_M(X^\prime,Y^\prime,Y^\prime,X^\prime)\leq R_N(X,Y,Y,X).
$$
  
The projection map $\tau:M\to M/G$ is a Riemannian submersion.
Therefore, we obtain the following fact: {\sl if a Lie group~$G$ acts
freely by isometries and the metric of~$M$ has positive sectional
curvature on all 2-planes orthogonal to the orbits of~$G$, then the
quotient metric on~$M/G$ has positive sectional curvature}.

In this article we prove a more general version of the above statement: 

\begin{theorem}                          
Let $M$ be a smooth manifold with an action of a group~$G$, the latter
being either~$S^1$ or~$S^3=Sp(1)$. Denote by~$N$ the fixed point set
of the action and by $\tau:M\to M/G$ the projection map. Let $W\subset
M$ be a neighborhood of~$N$. Assume that $M$ is endowed with a metric
$(\boldsymbol{\cdot}\,,\boldsymbol{\cdot})_0$ preserved by the action
of $G$, and the following conditions hold:
\begin{enumerate}
\item The metric $(\boldsymbol{\cdot}\,,\boldsymbol{\cdot})_0$ has positive 
sectional curvature on all 2-planes orthogonal to the orbits of $G$.
\item For each point $p\in M$, its isotropy group is either trivial or the 
whole $G$.
\item The space $N$ is a compact submanifold of~$M$ of codimension
$2\cdot(1+\dim G)$.
\end{enumerate}

Then the quotient space~$M/G$ posesses a smooth structure such
that~$\tau$ is smooth at the free points of the action; moreover, this
space carries a metric of positive sectional curvature that coincides
with the quiotient metric
of~$(\boldsymbol{\cdot}\,,\boldsymbol{\cdot})_0$ outside~$W$.
\end{theorem}

To construct the smooth structure on~$M/G$, we use infinite cones over
the Hopf fibration~$S^3/S^1\simeq S^2$ and the
fibration~$S^7/S^3\simeq S^4$. Some problems arise from the fact that
the quotient metric for the standard metric on the sphere in both of
the fibrations is the standard metric on the sphere of radius~$1/2$;
therefore, we have a conic singularity in the fixed point set $N$.

In the first section we give a construction which enables us to smooth
the quotient metric in a neighborhood of~$N$. In the second and third
sections we construct some particular metric on the normal bundle
over~$N$ and analyze its curvature. In the fourth section we show
that, under certain conditions, one can glue two metrics in a
neighborhood of a submanifold with positivity of sectional curvature
preserved. Finally, in the fifth section we construct the smooth
structure on~$M/G$ and the sought metric of positive sectional
curvature.

\section{Smoothing the Conic Metric}       

Consider the following metric on $\mathbb R^n$:
\begin{equation}\label{e:1.1}              
g=dr^2+g^2(r)d\varphi^2.
\end{equation}
Here $r(x)$ is the length of the vector~$x$, $d\varphi^2$ is the
standard metric of sectional curvature 1 on the sphere~$rS^{n-1}$, and
$g(r)$ is some positive function in~$C^2(0,\infty)$ with the third
derivative at zero.

The considered metric is $C^2$-smooth iff the following conditions hold:
\begin{equation}\label{e:1.2}              
g(0)=0,
\quad
g^\prime(0)=1,
\quad
g^{\prime\prime}(0)=0.
\end{equation}

By straightforward computations using the formulas for curvature of
Riemannian submersions, one can prove (this fact is, however, not used
in this paper directly) that the considered metric has positive
sectional curvature iff the following conditions hold:

\begin{equation}\label{e:1.3}              
g^{\prime\prime\prime}(0)<0,
\quad
g^{\prime\prime}(r)<0,
\quad
\big\vert g^\prime(r)
\big\vert<1\ \ \text{for}\ \ r>0.
\end{equation}

Now, we construct auxiliary functions~$g_\varepsilon$ to be used
later. Though the existence of such functions looks quite natural from
the geometrical point of view (see the figure below), the formal proof
faces certain difficulties; therefore, the proof is described at
length.

First, consider the function $g_0(r)=r(1-r^2)$ for $0\leq r<1/2$. It
satisfies~\eqref{e:1.2}, \eqref{e:1.3}, and the inequality
$g_0^\prime(r)\geq 1/4$. Next, for each~$\delta>0$ consider a {\it
mollifier}~$\omega_\delta(x)$, i.e., an even nonnegative smooth
function with support~$[-\delta,\delta]$ such that its integral over
the whole real axis equals~1. Given functions~$f$ and~$g$ defined on
the whole real axis, denote by~$f*g$ their convolution, i.e., the
integral
$$
(f*g)(x)=\int_{-\infty}^\infty f(t)g(x-t)\,dt.
$$

\begin{lemma}\label{l:1.1}\ \\           
\vskip-14pt\noindent
\begin{enumerate}
\item If a function $f$ is Lipshitz continuous and has the
derivative~$f^\prime$ almost everywhere, then
$(f*\omega_\delta)^\prime=f^\prime*\omega_\delta$.
\item If a function $f$ is concave, then such is the function
$f*\omega_\delta$.
\item $1*\omega_\delta=1$ and $x*\omega_\delta=x$.
\item If $f^{\prime\prime\prime}=0$ everywhere, then
$f*\omega_\delta=f+c$ for some constant~$c$.
\end{enumerate}
\end{lemma}

\begin{proof}
1. It suffices to apply the Lebesgue dominated convergence theorem.

2. The function $f$ is concave iff the following inequality holds for 
$x,y\in\mathbb R$ and $\lambda\in [0,1]$:
$$
f\big(\lambda x+(1-\lambda)y
 \big)\geq\lambda f(x)+(1-\lambda)f(y)
$$
If follows from the definition of convolution that the same inequality holds 
for the function~$f*w_\delta$.

3. The first statement is obvious. The second one follows from the formulas
$$
(x*\omega_\delta)(0)=0
\ \ \text{and}\ \ (x*\omega_\delta)^\prime=x^\prime*\omega_\delta.
$$

4. Since $f^\prime(x)$ is a linear function, we have 
$(f*\omega_\delta)^\prime(x)=f^\prime(x)$ and thus $f*\omega_\delta=f+c$.
\end{proof}

\begin{lemma}\label{l:1.2}           
For each $\varepsilon>0$, there exists a function $g_\varepsilon\in
C^2(0,1/2)$ with the following properties:
\begin{enumerate}
\item $g_\varepsilon$ satisfies~\eqref{e:1.2},
$g_\varepsilon^{\prime\prime}(r)\leq-r$, and
$g_\varepsilon^\prime(r)\leq g_0^\prime(r)$.
\item If $r\geq\varepsilon$, then $g_\varepsilon(r)=g_0(r)/2$.
\end{enumerate}
\end{lemma}

\begin{figure}\centering
\includegraphics{dyatlov.1}
\end{figure}

\begin{proof} We may assume that $\varepsilon<1/2$. Consider the
set~$\mathcal A$ of all functions of class $C^2[0,\varepsilon]$ that
satisfy the relations
\begin{gather}
g(0)=0,
\ \ g^{\prime\prime}(0)=0,
\ \ g^{\prime\prime}(r)\leq -r,\label{e:1.4}\\              
g^\prime(r)\leq g_0^\prime(r)+g^\prime(0)-1.\label{e:1.5}   
\end{gather}
It is clear that $\mathcal A$ is convex. Consider the operator
$T_\varepsilon:\mathcal A\to\mathbb R^4$ defined as follows:
$$
T_\varepsilon(g)=
\big(g(\varepsilon),
     g^\prime(\varepsilon),
     g^{\prime\prime}(\varepsilon),
     g^\prime(0)
\big).
$$
It suffices to prove that the set $T_\varepsilon(\mathcal A)$ contains
the point
$$
A=\left({g_0(\varepsilon)\over 2},
        {g_0^\prime(\varepsilon)\over 2},
        {g_0^{\prime\prime}(\varepsilon)\over 2},1
  \right).
$$
Since this set is convex and has finite dimension, it suffices to
check that it is dense in a neighborhood of the point~$A$. Consider
the auxiliary function
$\varphi(r)=(r^3-3\varepsilon^2r+2\varepsilon^3)/6$. Note that
$\varphi^{\prime\prime}(r)=r$ and
$\varphi(\varepsilon)=\varphi^\prime(\varepsilon)=\nobreak0$.

Take a point $B=(b_0,b_1,b_2,b_3)$. If~$B$ and~$A$ are sufficiently close to each other,
then the following conditions hold:
\begin{equation}\label{e:1.6}       
\begin{gathered}
b_2+\varphi^{\prime\prime}(\varepsilon)<0,
\ \ \varepsilon b_1+\varphi(0)<b_0<(b_3-1)\varepsilon+g_0(\varepsilon),\\
b_1<b_3-1+g_0^\prime(\varepsilon).
\end{gathered}
\end{equation}

In order to prove that the set $T_\varepsilon(\mathcal A)$ is dense in
a neighborhood of the point~$A$, we construct a one-parameter family
of functions $g_\delta\in\mathcal A$ such that
$T_\varepsilon(g_\delta)\to B$ when $\delta\to +0$. Construct a
function $g_1\in C^2[0,\varepsilon]$ with the properties

\begin{equation}\label{e:1.7}          
\begin{gathered}
g_1^{\prime\prime}\leq 0,
\ \ g_1(\varepsilon)=b_0,
\ \ g_1^\prime(\varepsilon)=b_1,
\ \ g_1^{\prime\prime}(\varepsilon)=b_2+\varphi^{\prime\prime}(\varepsilon),\\
g_1(0)>\varphi(0),
\ \ g_1^\prime(r)-\varphi^\prime(r)\leq b_3-1+g_0^\prime(r).
\end{gathered}
\end{equation}

It can be sought in the form
$$
g_1(r)=b_0+b_1(r-\varepsilon)-\int_0^{\varepsilon-r} u(t)\,dt.
$$
The function $u\in C^1[0,\varepsilon]$ must satisfy the following
conditions:
$$
\gathered
u(0)=0,
\ \ u^\prime(0)=-\big(b_2+\varphi^{\prime\prime}(\varepsilon)
                 \big),
\ \ u^\prime\geq 0,\\
\int_0^\varepsilon u(t)\,dt<b_0-b_1\varepsilon-\varphi(0),
\ \ u(\varepsilon-r)\leq\varphi^\prime(r)+g_0^\prime(r)+b_3-1-b_1.
\endgathered
$$
Such function $u$ exists if the conditions~\eqref{e:1.6} hold.

We have
$g_1(\varepsilon)-\varphi(\varepsilon)<(b_3-1)\varepsilon+g_0(\varepsilon)$
and at the same time $g_1(0)-\varphi(0)>g_0(0)$. Therefore, there is a
point $r_0\in(0,\varepsilon)$ such that
$$
g_1(r_0)-\varphi(r_0)=(b_3-1)r_0+g_0(r_0).
$$
Consider the following function~$g_2$ defined on the whole real axis:
$$
g_2(r)=
\begin{cases}
(b_3-1)r+g_0(r),  &r<r_0,\\
g_1(r)-\varphi(r),&r_0\leq r<\varepsilon,\\
b_0+b_1(r-\varepsilon)+\dfrac{b_2}{2}(r-\varepsilon)^2,
                  &r\geq\varepsilon.
\end{cases}
$$
It is obvious that this function is continuous everywhere
and~$C^2$-smooth at all points but~$r_0$. Take the convolution
$g_\delta=g_2*\omega_\delta$. It is well-known that $g_\delta\in
C^\infty$, and $T_\varepsilon(g_\delta)\to T_\varepsilon(g_2)=B$ as
$\delta\to 0$. It remains to prove that $g_\delta\in\mathcal A$ for
small~$\delta$. The equalities
$g_\delta(0)=g_\delta^{\prime\prime}(0)=0$ follow from
Lemma~\ref{l:1.1}~(item 1) and from the fact that $g_2$ is odd in a
neighborhood of zero. The inequality $g_\delta^{\prime\prime}(r)\leq
-r$ holds in a neighborhood of zero due to Lemma~\ref{l:1.1}~(items 1
and~3), and outside this neighborhood this inequality follows from
Lemma~\ref{l:1.1}~(item 2) and the concavity of the
function~$g_2+\varphi$ for~$r>0$.

We now show that $g_\delta$ satisfies~\eqref{e:1.5}. In a neighborhood
of the segment~$[0,\varepsilon]$ we have $g_2^\prime\leq
b_3-1+g_0^\prime$, so Lemma~\ref{l:1.1}~(items 1 and~4) yields
\begin{equation*}
g_\delta^\prime=
g_2^\prime*\omega_\delta\leq (b_3-1+g_0^\prime)*\omega_\delta=
g_0^\prime-1+
g_\delta^\prime(0).
\end{equation*}
\end{proof}

\begin{lemma}\label{l:1.3}           
Let $0<t<1/3$ and $g_0(t_0)=g_\varepsilon(t)$. Then
$g_\varepsilon^\prime(t)\leq g_0^\prime(t_0)$ and
$7tg_\varepsilon^\prime(t)\geq t_0g_0^\prime(t_0)+g_\varepsilon(t)$.
\end{lemma}

\begin{proof} We have $g_\varepsilon^\prime\leq g_0^\prime$;
therefore, $g_\varepsilon\leq g_0$ and $t_0\leq t$. Thus
$g_\varepsilon^\prime(t)\leq g_0^\prime(t)\leq g_0^\prime(t_0)$.

The functions $g_0$ and $g_\varepsilon$ are concave; moreover,
$g_0^\prime(0)=g_\varepsilon^\prime(0)=1$, $g_0(1/2)=3/8$, and
$g_\varepsilon(1/2)=3/16$; therefore, for $0\leq t\leq 1/2$ we have
$$
{3\over 4}t\leq g_0(t)\leq t,
\quad
{3\over 8}t\leq g_\varepsilon(t)\leq t.
$$
Now,
\begin{equation*}
t_0g_0^\prime(t_0)\leq
t_0\leq
{4\over 3}g_0(t_0)=
{4\over 3}g_\varepsilon(t)\leq
{4\over 3}t,
\quad
tg_\varepsilon^\prime(t)\geq
tg_\varepsilon^\prime(1/3)\geq
{t\over 3}.
\end{equation*}
\end{proof}

\section{Certain Metrics on Vector Bundles}   

Some results of these section follow from~\cite[theorem~3.5]{Vi} and
\cite[theorem~9.59]{Be}. However, we state them together with proofs
for completeness of exposition.

Suppose that $N$ is an $n$-dimensional compact manifold and $M$ is
vector bundle over~$N$ of rank~$m-n$; i.e.,
\begin{equation}\label{e:2.1}            
M=\big\{(p,v)\mid p\in N,\ v\in V_p
  \big\}.
\end{equation}
Here $V_p$ is an ($m-n$)-dimensional vector space, which we call the
{\it fiber\/} at the point~$p$. We assume that each~$V_p$ is endowed
with an inner product
$\langle\boldsymbol{\cdot}\,,\boldsymbol{\cdot}\rangle$ depending
smoothly on $p$. We can embed~$N$ into~$M$ as~$\{v=\nobreak0\}$.
Denote by $\pi:M\to N$ the projection map. A line of the
form~$(p,tv)_{t\geq 0}$ is called a {\it ray emanating from~$p$ in the
direction~$v$}.

Define the function $r:M\to\mathbb R$ as follows:
\begin{equation}\label{e:2.2}            
r(p,v)=\sqrt{\langle v,v\rangle}.
\end{equation}

It is obvious that the function~$r$ is continuos and belongs
to~$C^\infty(M^\circ)$, where $M^\circ=M\setminus N$ is the regular
part of~$M$. A {\it neighborhood of}~$N$ is a submanifold
$\{r<\rho\}\subset M$ for~$\rho>0$.

Note that if $K=\{r=1\}\in M$, then there is a standard identification
$M^\circ\simeq K\times (0,\infty)$. That makes it possible to identify
different {\it radial fibers}, i.e., the submanifolds of the kind
$\{r=t\}$ for $t>0$. Denote by~$\partial/\partial r$ the image of the
standard vector field~$1$ on~$(0,\infty)$ under this identification.
The images of vector fields on~$K$ under this identification are
called {\it radial basic fields\/}. {\it Vertical fields\/} are the radial
basic fields that are tangent to the fibers. Note that for each radial
basic field~$A$ we have

\begin{equation}\label{e:2.3}              
\bigg[\pop r,A
\bigg]=0.
\end{equation}

Moreover, the Lie bracket of each two radial basic fields is a 
radial basic field.

A {\it trivialization\/} of the bundle~$M$ is a pair~$(U,\Theta)$,
where $U$ is an open subset in~$N$ and $\Theta$ is a diffeomorphism of
$U\times\mathbb R^{m-n}$ onto~$\pi^{-1}(U)$ such that for each~$p\in
U$, the restriction of~$\Theta$ to $p\times\mathbb R^{m-n}$ is a
linear isomorphism onto~$V_p$ that maps the standard inner product
on~$\mathbb R^{m-n}$ to the inner
product~$\langle\boldsymbol{\cdot}\,,\boldsymbol{\cdot}\rangle$
on~$V_p$. We~assume that the bundle~$M$ is {\it locally trivial};
i.e., for each~$p\in N$, there exists a trivialization~$(U,\Theta)$
such that~$p\in U$.

Let~$(U_1,\Theta_1)$ and~$(U_2,\Theta_2)$ be two trivializations.
Then there is a smooth map~$\Psi_{12}:U_1\cap U_2\to SO(m-n)$
such that for $p\in U_1\cap U_2$ and $x\in\mathbb R^{m-n}$ we have

\begin{equation}\label{e:2.4}                
\Theta_2(p,x)=\Theta_1(p,\Psi_{12}(p)\cdot x).
\end{equation}

Consider a metric of the form $dr^2+g^2(r)d\varphi^2$ on $\mathbb
R^{m-n}$ as in~\eqref{e:1.1}. Note that this metric is invariant under
the standard action of the group~$SO(m-n)$. We call a vector field $V$
on~$\mathbb R^{n-m}$ a {\it turn field\/} iff it has the
form~$V(x)=A\cdot x$ for some matrix $A\in\mathfrak{so}(m-n)$. Note
that turn fields are Killing fields of the action of~$SO(m-n)$, so we
have the following relation for vector fields~$A$ and~$B$ and a turn
field~$C$:
\begin{equation}\label{e:2.5}              
(\nabla_AC,B)+(\nabla_BC,A)=0.
\end{equation}
Furthermore, a turn field is mapped by each element of~$SO(m-n)$ to
another turn field, so~\eqref{e:2.4} makes it possible to define turn
fields on each fiber~$V_p$. A vector field on~$M$ is called a turn
field iff it is tangent to the fibers and is a turn field on each
fiber.

Consider a trivialization~$(U,\Theta)$. We call the images of vector
fields on~$U$ under~$\Theta$ {\it locally horizontal}. The images
under~$\Theta$ of vector fields on~$\mathbb R^{m-n}$ that are constant
in the Cartesian coordinates are called {\it locally vertical\/}
fields on~$M$. For a vector field $X$ on~$M$, denote by~$X^\prime$ its
orthogonal projection to the fibers with respect to some metric
on~$M$.

\begin{lemma}\label{l:2.1}              
Suppose that $(U,\Theta_1)$ and $(U,\Theta_2)$ are two trivializations
and $X_1$ and $X_2$ are the locally horizontal fields in the first and
second trivializations, respectively, that are produced by the same
vector field on~$U$. Then~$X_1-X_2$ is a turn field.
\end{lemma}

\begin{proof}
Consider the field $X$ on~$U$ that produces both~$X_1$ and~$X_2$. Let
$\Psi_{12}$ be the map that connects the two trivializations
by~\eqref{e:2.4}. Then
$$
(X_2-X_1)_{\Theta_1(p,x)}=
d\Theta_{1(p,x)}\big(0,(X\Psi_{12})(p)\cdot\Psi_{12}(p)^{-1}\cdot x
                \big).
$$
Since $\Psi_{12}(p)\in SO(m-n)$, we have
$(X\Psi_{12})(p)\cdot\Psi_{12}(p)^{-1}\in\mathfrak{so}(m-n)$ and the
lemma follows.
\end{proof}

This lemma shows that the following definition is correct:

\begin{definition}\label{d:2.1}           
Let $g$ be a function satisfying~\eqref{e:1.2}. A metric
$(\boldsymbol{\cdot}\,,\boldsymbol{\cdot})$ on~$M$ is called
{\it proper\/} or {\it $g$-proper} iff the following conditions hold:
\begin{enumerate}
\item The restriction of the metric
$(\boldsymbol{\cdot}\,,\boldsymbol{\cdot})$ to each fiber~$V_p$ is
equal to~$dr^2+g^2(r)d\varphi^2$.
\item The orthogonal projection to the fibers of each vector field that
is locally horizontal in some trivialization is a turn field.
\end{enumerate}
\end{definition}

\begin{lemma}\label{l:2.2}                  
Let $(\boldsymbol{\cdot}\,,\boldsymbol{\cdot})$ be a proper metric.
Then:
\begin{enumerate}
\item At every point~$p\in N$, the spaces~$T_pN$ and~$T_pV_p$ are orthogonal.
\item Each fiber~$V_p$ is totally geodesic.
\end{enumerate}
\end{lemma}

\begin{proof}
1. At the points of~$N$, the orthogonal projection of every locally
horizontal field to the fiber is zero.

2. Consider a trivialization~$(U,\Theta)$. Suppose that~$A$ and $B$
are two locally vertical vector fields and $X$ is a locally horizontal
field. It suffices to prove that $(\nabla_AB,X-X^\prime)=0$. We have
\begin{gather*}
[A,X]=[B,X]=[A,B]=0,
\quad X(A,B)=0,\\
(\nabla_AX^\prime,B)+(\nabla_BX^\prime,A)=0,\\
(A,X)=(A,X^\prime),
\quad (B,X)=(B,X^\prime),\\
2(\nabla_AB,X)=A(B,X)+B(A,X)-X(A,B),\\
2(\nabla_AB,X^\prime)=
  A(B,X^\prime)+B(A,X^\prime)-(\nabla_AX^\prime,B)-(\nabla_BX^\prime,A).
\end{gather*}
The proof is complete.
\end{proof}

Consider a trivialization~$(U,\Theta)$. Assume that we are given:
\begin{enumerate}
\item A linear map $Q_p:T_pN\to\mathfrak{so}(m-n)$ that depends smoothly on the 
point~$p\in N$;

\item A function~$g$ satisfying~\eqref{e:1.2};

\item A metric~$g_N$ on~$N$.
\end{enumerate}

\begin{lemma}\label{l:2.3}                  
There exists a unique $g$-proper metric
$(\boldsymbol{\cdot}\,,\boldsymbol{\cdot})$ on $\pi^{-1}(U)$ such that
\begin{enumerate}
\item For $p\in U$, $x\in\mathbb R^{m-n}$, and a locally horizontal 
field~$X$ we have
\begin{equation}\label{e:2.6}        
X^\prime_{\Theta(p,x)}=d\Theta_{(p,x)}\big(0,Q_p(X_p)\cdot x
                                      \big).
\end{equation}
\item The map $\pi$ is a Riemannian submersion onto~$N$ endowed with the 
metric~$g_N$.
\end{enumerate}

Furthermore, the submanifold~$N$ is totally geodesic with respect to
this metric, and the following relations hold at any point~$p\in N$
for all locally vertical fields~$A$, $B$, and $C$ and locally
horizontal fields~$X$ and $Y$:

\begin{equation}\label{e:2.7}            
A(B,C)=0,
\ \ A(X,Y)=0,
\ \ A(B,X)=\big\langle Q_p(X_p)\cdot A,B
           \big\rangle.
\end{equation}

The fields~$A$ and~$B$ on the right-hand side of the last equation are
identified with the vector fields on~$\mathbb R^{m-n}$ that produces
them.

Moreover, if we identify a radial fiber~$\{r=t\}$ with~$K=\{r=1\}$ in
a standard way, then the
metric~$(\boldsymbol{\cdot}\,,\boldsymbol{\cdot})$ depends only on the
map~$Q_p$, the metric~$g_N$, and the value~$g(t)$ (and does not
explicitly depend on~$t$).
\end{lemma}

\begin{proof} The fact that the metric is $g$-proper determines its
form on the fibers~$V_p$. The
metric~$(\boldsymbol{\cdot}\,,\boldsymbol{\cdot})$ is given by
formulas
\begin{equation}\label{e:2.8}           
(A,X)=(A,X^\prime),
\quad (X,Y)=(X^\prime,Y^\prime)+\pi^*g_N(X,Y).
\end{equation}
Here the field~$A$ is locally vertical, the fields~$X$ and~$Y$ are
locally horizontal, and the fields~$X^\prime$ and $Y^\prime$ are
determined from~\eqref{e:2.6}. Since $Q_p$ is a linear map from~$T_pN$
to~$\mathfrak{so}(m-n)$, \eqref{e:2.8} defines a smooth proper
Riemannian metric.

The equalities~\eqref{e:2.7} can be checked directly. The fact
that~$N$ is totally geodesic is equivalent to the following equation
holding at the points
of~$N$:
\begin{equation}\label{e:2.9}        
(\nabla_XX,A)=0.
\end{equation}
But at these points we have
\begin{equation}\label{e:2.10}       
(\nabla_XX,A)=X(A,X)-
(X,\nabla_XA)=-
(X,\nabla_AX)=-{1\over 2}A(X,X).
\end{equation}
It follows that~\eqref{e:2.9} is equivalent to the second equation
in~\eqref{e:2.7}. The last statement of the lemma is verified directly
using \eqref{e:2.8} and the form of the metric on the fibers.
\end{proof}

Let $M_0$ be a manifold endowed with a Riemannian
metric~$(\boldsymbol{\cdot}\,,\boldsymbol{\cdot})_0$; let $N$ be its
compact totally geodesic submanifold. Define~$g_N$ as the restriction
of the metric~$(\boldsymbol{\cdot}\,,\boldsymbol{\cdot})_0$ to~$N$.
Consider as~$M$ the {\it normal bundle\/} over~$N$; i.e., let
$$
V_p=T_pN^\perp=\{v\in T_pM_0\mid v\perp T_pN\}.
$$

Take the inner product
$\langle\boldsymbol{\cdot}\,,\boldsymbol{\cdot}\rangle$ on the
fiber~$V_p$ to be the restriction of the metric
$(\boldsymbol{\cdot}\,,\boldsymbol{\cdot})_0$ to~$T_pN^\perp$.
Consider the map~$\Phi:M\to M_0$ defined as follows:
\begin{equation}\label{e:2.11}       
\Phi(p,v)=\exp_pv.
\end{equation}
It is well-known that this map is a diffeomorphism for small~$r$.
Using~$\Phi$, we pull back the
metric~$(\boldsymbol{\cdot}\,,\boldsymbol{\cdot})_0$ to a neighborhood
of~$N$ in~$M$.

\begin{lemma}\label{l:2.4}              
Let $g$ be a function satisfying~\eqref{e:1.2}. Then there exists
unique $g$-proper metric~$(\boldsymbol{\cdot}\,,\boldsymbol{\cdot})_1$
on~$M$ such that the 1-jets of the
metrics~$(\boldsymbol{\cdot}\,,\boldsymbol{\cdot})_0$
and~$(\boldsymbol{\cdot}\,,\boldsymbol{\cdot})_1$ coincide at the
points of $N$, and $\pi$ is a Riemannian submersion, if we take the
metric~$g_N$ on~$N$. Furthermore, if we identify some radial
fiber~$\{r=t\}$ with~$K$, then the
metric~$(\boldsymbol{\cdot}\,,\boldsymbol{\cdot})_1$ on this fiber
depends only on the
metric~$(\boldsymbol{\cdot}\,,\boldsymbol{\cdot})_0$ and the
value~$g(t)$ (and does not depend explicitly on~$t$). Moreover, for
each field~$X$, its orthogonal projection to the fibers~$X^\prime$ is
independent of the function~$g$.
\end{lemma}

\begin{proof} It suffices to prove local existence and uniqueness of
the metric~$(\boldsymbol{\cdot}\,,\boldsymbol{\cdot})_1$. Let
$(U,\Theta)$ be some trivialization on~$M$. Suppose that the vector
fields~$A$, $B$, and $C$ are locally vertical and~$X$ and~$Y$ are
locally horizontal. It is clear that, if the sought metric exists,
then it may be constructed using Lemma~\ref{l:2.3} for some map~$Q_p$.
In this case the restrictions of the
metrics~$(\boldsymbol{\cdot}\,,\boldsymbol{\cdot})_0$
and~$(\boldsymbol{\cdot}\,,\boldsymbol{\cdot})_1$ to~$N$ conincide
automatically, so in order to make the $1$-jets of the two metrics
coincide at the points of $N$, it suffices to obtain the
equalities~\eqref{e:2.7} together with the following equalities for
the metric~$(\boldsymbol{\cdot}\,,\boldsymbol{\cdot})_0$ at the points
of~$N$:
\begin{equation}\label{e:2.12}         
(A,X)_0=0,\quad
(A,B)_0=\langle A,B\rangle
\end{equation}
It can be seen from the construction of the normal bundle that~\eqref{e:2.12} 
and the first equation in~\eqref{e:2.7} hold. The second equation in~\eqref{e:2.7}
follows from the fact that~$N$ is totally geodesic
(see~\eqref{e:2.9} and \eqref{e:2.10}). The last equation in~\eqref{e:2.7}
is equivalent to the following one:
\begin{equation}\label{e:2.13}      
\big\langle Q_p(X_p)\cdot A, B
\big\rangle=A(B,X)_0.
\end{equation}
But the expression on the right-hand side depends multilinearly
on~$A$, $B$, and~$X$, so~\eqref{e:2.13} defines the unique map~$Q_p$. It
remains to prove that the map~$Q_p$ defined above takes its values
in~$\mathfrak{so}(m-n)$, i.e., that the following equality holds at
the points of~$N$:
$$
A(A,X)_0=0.
$$
Indeed, $\nabla_AA=0$ at the points of~$N$, so
$$
A(A,X)_0=(A,\nabla_AX)_0=(A,\nabla_XA)_0={1\over 2}X(A,A)_0=0.
$$
The last two statements of this lemma follow from the last statement of the 
previous lemma, \eqref{e:2.6} and~\eqref{e:2.13}.
\end{proof}

We now establish some properties of proper metrics that will be useful
in computing the curvature. Consider a vector field~$X$ on~$N$. There
exists a unique field~${}\,\widehat{\!X}$ on~$M$ which is $\pi$
related to~$X$ and orthogonal to the fibers; we call it a {\it basic
field} (produced by the field~$X$). Note that the
vector~$\partial/\partial r$ is orthogonal to the radial fibers, so
all basic fields are tangent to the radial fibers.

\begin{lemma}\label{l:2.5}           
If~$(\boldsymbol{\cdot}\,,\boldsymbol{\cdot})$ is a proper metric on~$M$, then
\begin{enumerate}
\item For all radial basic fields~$A$ and~$X$ such that~$A$
is tangent to the fibers, the value~$(A,X)/(A,A)$ is constant on every ray.
\item If a field~$X$ is basic, then it is radial basic.
\end{enumerate}
\end{lemma}

\begin{proof}
1. Consider a trivialization~$(U,\Theta)$. For radial basic fields~$X$
which are tangent to the fiber, the statement follows from the form
of the metric on the fibers, so we may assume that~$X$ is locally
horizontal. Then $(A,X)=(A,X^\prime)$ and~$X^\prime$ is a turn field;
therefore, it is tangent to the fibers and is radial basic.

2. Consider the radial basic field~$Y$ that coincides with~$X$ on the
radial fiber $\{r=t\}$ for some $t$. This is possible since~$X$ is
tangent to the radial fibers. Now,~$X$ and~$Y$ are $\pi$~related to
the same field on~$N$, and~$Y$ is orthogonal to the fibers due to the
first statement of this lemma. Therefore, $X=Y$.
\end{proof}

\begin{lemma}\label{l:2.6}         
Let $(\boldsymbol{\cdot}\,,\boldsymbol{\cdot})$ be a proper metric and
let $A$ and $B$ be radial basic fields. Denote by~$\II$ the second
fundamental form of the radial fibers. Then the following equalities
hold:
\begin{gather}
\bigg(\II(A,B),\pop r
\bigg)=-{1\over 2}\pop r(A,B),\label{e:2.14}\\  
\bigg(\DoP rA,B
\bigg)={1\over 2}\pop r(A,B),\label{e:2.15}\\   
\bigg(\DoP rA,\pop r
\bigg)=0.\label{e:2.16}                         
\end{gather}
\end{lemma}

\begin{proof}
The statement of the lemma follows from the equalities
\begin{gather*}
\bigg(\II(A,B),\pop r
\bigg)=\bigg(\nabla_AB,\pop r
       \bigg)
=\bigg(\nabla_BA,\pop r
 \bigg),\\
 \bigg(\nabla_AB,\pop r
 \bigg)=A\bigg(B,\pop r
         \bigg)-\bigg(B,\nabla_A\pop r
                \bigg)
=-\bigg(B,\DoP rA
  \bigg),\\
\pop r(A,B)=\bigg(A,\DoP rB
            \bigg)+\bigg(B,\DoP rA
                   \bigg),\\
\bigg(\DoP rA,\pop r
\bigg)=\bigg(\nabla_A\pop r,\pop r
       \bigg)={1\over 2}A
       \bigg(\pop r,\pop r
       \bigg)=0.
\end{gather*}
\end{proof}

\section{The Curvature of a Certain Family of Proper Metrics}  

As in Lemma~\ref{l:2.4}, suppose that $M_0$ is a manifold with a given
metric~$(\boldsymbol{\cdot}\,,\boldsymbol{\cdot})_0$, $N$ is its
totally geodesic submanifold, and $M$ is the normal bundle over $N$.
The metric~$(\boldsymbol{\cdot}\,,\boldsymbol{\cdot})_0$ is pulled
back to~$M$ via the map~$\Phi$ from~\eqref{e:2.11}. We also assume
that~$g$ is a given function satisfying \eqref{e:1.2}, and $L\geq 0$
is a constant. Using Lemma~\ref{l:2.4}, consider the $g$ proper
metric~$(\boldsymbol{\cdot}\,,\boldsymbol{\cdot})_1$ on~$M$ such that
the map $\pi$ is a Riemannian submersion. Construct the
metric~$(\boldsymbol{\cdot}\,,\boldsymbol{\cdot})$ in a neighborhood
of~$N$ by the formula
\begin{equation}\label{e:3.1}                      
(X,Y)=(X,Y)_1-Lr^2\cdot(\pi^*g_N)(X,Y).
\end{equation}
Here $X$ and~$Y$ are arbitrary vectors. It is clear that the constructed
metric is smooth. It is positive definite for $r^2<1/L$ and its
$1$-jets at the points of~$N$ coincide with those of the
metric~$(\boldsymbol{\cdot}\,,\boldsymbol{\cdot})_0$. Furthermore, on
each radial fiber~$\{r=t\}$ the
metric~$(\boldsymbol{\cdot}\,,\boldsymbol{\cdot})$ is determined by
the metric~$(\boldsymbol{\cdot}\,,\boldsymbol{\cdot})_0$, the
constant~$L$, and the value~$g(t)$. Moreover, for each field~$X$ its
orthogonal projection~$X^\prime$ to the fibers depends only on the
original metric~$(\boldsymbol{\cdot}\,,\boldsymbol{\cdot})_0$ (that
is, it does not depend on the function~$g$ and the constant~$L$).
Thus, the metric~$(\boldsymbol{\cdot}\,,\boldsymbol{\cdot})_0$ defines
what fields are basic in the
metric~$(\boldsymbol{\cdot}\,,\boldsymbol{\cdot})$.

The map~$\pi$ restricted to each radial fiber is a Riemannian
submersion with totally geodesic fibers from~$M$ endowed with the
metric~$(\boldsymbol{\cdot}\,,\boldsymbol{\cdot})$ to~$N$ endowed with
the metric~$(1-Lr^2)g_N$. Therefore (see~\cite[p.~9\,C]{Be}), for all
basic fields~$X$ and $Y$ and vertical fields~$A$ and $B$ we have
\begin{equation}
\begin{gathered} \label{e:3.2}                    
[A,X]\in T(V_p),
\quad
(\nabla_XY,A)={1\over 2}\big(A,[X,Y]
                        \big),\\
(\nabla_AX,B)=0,
\quad
(\nabla_XA,B)=\big([X,A],B
              \big),\\
(\nabla_XA,Y)=-{1\over 2}\big(A,[X,Y]
                         \big)=(\nabla_AX,Y).
\end{gathered}
\end{equation}

Denote the covariant derivative and curvature tensor of the radial fibers by 
the symbols~$\overline\nabla$ and~$\olr$. It is well-known that if the
fields~$A$ and~$B$ are tangent to the radial fibers, then

\begin{gather}\label{e:3.3}                            
\nabla_AB=\overline\nabla_AB+\II(A,B),\\
\begin{align}
R(A,B,B,A)=\olr(A,B,B,A)\nonumber\\
          \qquad-\big(\II(A,A),\II(B,B)
                  \big)+\big(\II(A,B),\II(A,B)
                        \big).
\end{align}\label{e:3.4}                             
\end{gather}

Suppose that the fields $A$ and~$B$ are vertical, and the fields~$X$ and
$Y$ are basic. It follows from~\eqref{e:2.15} that

\begin{gather*}
\bigg(\DoP rA,B
\bigg)={1\over 2}\pop r(A,B)={g^\prime(r)\over g(r)}(A,B),\\
\bigg(\DoP rX,A
\bigg)=\bigg(\DoP rA,X
       \bigg)={1\over 2}\pop r(X,A)=0,\\
\bigg(\DoP rX,Y
\bigg)={1\over 2}\pop r(X,Y)=-{Lr\over 1-Lr^2}(X,Y),
\end{gather*}

From here and~\eqref{e:2.16} we obtain

\begin{gather}
\DoP rA={g^\prime(r)\over g(r)}A,\label{e:3.5}\\     
\DoP rX=-{Lr\over 1-Lr^2}X.      \label{e:3.6}       
\end{gather}

We are now ready to compute the sectional curvature of the
metric~$(\boldsymbol{\cdot}\,,\boldsymbol{\cdot})$. To this end,
consider two vector fields
$$
E=a\pop r+A+X,
\quad
F=b\pop r+B+Y.
$$

We assume that the vector fields $A$ and $B$ are vertical, the fields
$X$ and $Y$ are basic, and $[A,B]=0$. We also assume that the
vectors~$E$ and $F$ form an orthonormal system at the point where the
curvature is computed. We introduce the fields
\begin{equation}
\begin{gathered}
P=A+X,
\quad Q=B+Y,
\quad\widehat{A}=a\pop r+A,
\quad\widehat{B}=b\pop r+B,\\
C=aB-bA,
\quad Z=aY-bX,
\quad G=C+Z.
\end{gathered}
\end{equation}

We have

\begin{equation}\label{e:3.7}            
R(E,F,F,E)=R(P,Q,Q,P)
+R\bigg(\pop r,G,G,\pop r
  \bigg)-2R\bigg(\pop r,G,P,Q
           \bigg).
\end{equation}

For each radial basic field~$H$ we have

\begin{align*}
R\bigg(\pop r,H,H,\pop r\bigg)
&=\bigg({D\over\partial r}\nabla_HH,\pop r
  \bigg)-\bigg(\nabla_H\DoP rH,\pop r
         \bigg)\\
&=\pop r\bigg(\nabla_HH,\pop r
        \bigg)+\bigg\vert\DoP rH
               \bigg\vert^2\\
&=\bigg\vert\DoP rH
  \bigg\vert^2-{1\over 2}{\partial^2\over\partial r^2}(H,H).
\end{align*}

From here, \eqref{e:3.5}, and~\eqref{e:3.6} we obtain

\begin{equation}\label{e:3.8}                   
\begin{gathered}
R\bigg(\pop r,G,G,\pop r
 \bigg)=-{g^{\prime\prime}(r)\over g(r)}(C,C)+
 \bigg({L\over 1-Lr^2}+{(Lr)^2\over(1-Lr^2)^2}
 \bigg)(Z,Z).
\end{gathered}
\end{equation}

Furthermore,

\allowdisplaybreaks
\begin{align*}
R\bigg(\!A,B,C,\pop r\!
 \bigg)
&\,{=}\,A\bigg(\!\nabla_BC,\pop r\!
         \bigg){-}B
         \bigg(\!\nabla_AC,\pop r\!
         \bigg){+}
         \bigg(\!\nabla_AC,\DoP rB\!
         \bigg){-}
         \bigg(\!\nabla_BC,\DoP rA\!
         \bigg)\\
&={1\over 2}\pop r\big(B(A,C){-}A(B,C)
                  \big){+}
                  \bigg(\!\nabla_AC,\DoP rB\!
                  \bigg){-}
                  \bigg(\!\nabla_BC,\DoP rA\!
                  \bigg)\\
&={g^\prime(r)\over g(r)}\big(B(A,C)-A(B,C)+(\nabla_AC,B)-(\nabla_BC,A)
                         \big)\\
&=0,\\
R\bigg(\pop r,Z,X,Y\!
 \bigg)
&=\pop r(\nabla_ZX,Y){-}Z
                 \bigg(\!\DoP rX,Y\!
                 \bigg){+}
                 \bigg(\!\DoP rX,\nabla_ZY\!
                 \bigg){-}
                 \bigg(\!\DoP rY,\nabla_ZX\!
                 \bigg)\\
&={Lr\over 1{-}Lr^2}
                 \big(-2(\nabla_ZX,Y){+}Z(X,Y){-}(\nabla_ZY,X){+}(\nabla_ZX,Y)
                 \big)\\
&=0,\\
R\bigg(A,B,Z,\pop r
 \bigg)
&={1\over 2}\pop r\big(B(A,Z){-}A(B,Z)
                  \big){+}
                 \bigg(\!\nabla_AZ,\DoP rB\!
                 \bigg){-}
                 \bigg(\!\nabla_BZ,\DoP rA\!
                 \bigg)\\
&={g^\prime(r)\over g(r)}
                 \big((\nabla_AZ,B)-(\nabla_BZ,A)
                 \big)=0,\\
R\bigg(\pop r,C,A,Y
 \bigg)
&=\bigg({D\over\partial r}\nabla_CA,Y
         \bigg)-
         \bigg(\nabla_C\DoP rA,Y
  \bigg)=0,\\
R\bigg(\pop r,C,X,Y
 \bigg)
&=\pop r(\nabla_CX,Y){-}C
                 \bigg(\!\DoP rX,Y\!
                 \bigg){+}
                 \bigg(\!\DoP rX,\nabla_CY\!
                 \bigg){-}
                 \bigg(\!\DoP rY,\nabla_CX\!
                 \bigg)\\
&=-\bigg({Lr\over 1-Lr^2}+{g^\prime(r)\over g(r)}
         \bigg)\big(C,[X,Y]
               \big),\\
R\bigg(\pop r,Z,A,Y\!
  \bigg)
&=\pop r(\nabla_ZA,Y){-}Z
         \bigg(\DoP r A,Y\!
         \bigg)
{+}\bigg(\DoP rA,\nabla_ZY\!
   \bigg){-}\bigg(\DoP rY,\nabla_ZA
            \bigg)\\
&={1\over 2}\pop r\big(A,[Y,Z]
                        \big)+{1\over 2}
        \bigg({Lr\over 1-Lr^2}-{g^\prime(r)\over g(r)}
        \bigg)
        \big(A,[Y,Z]
        \big)\\
&={1\over 2}
        \bigg({Lr\over 1-Lr^2}+{g^\prime(r)\over g(r)}
        \bigg)\big(A,[Y,Z]
              \big).
\end{align*}

Hence, we find that

\begin{align}\label{e:3.9}         
&R\bigg(\pop r,G,P,Q
  \bigg)\nonumber\\
&\qquad={1\over 2}\bigg({Lr\over 1-Lr^2}+{g^\prime(r)\over g(r)}
                  \bigg)\Big(-2\big(C,[X,Y]
                               \big)+\big(A,[Y,Z]
                                     \big)-\big(B,[X,Z]
                                           \big)
                        \Big)\nonumber\\
&\qquad=-{3\over 2}\bigg({Lr\over 1-Lr^2}+{g^\prime(r)\over g(r)}
                   \bigg)\big(C,[X,Y]
                         \big).
\end{align}

Finally, from~\eqref{e:2.14} and~\eqref{e:3.4} we conclude that

\allowdisplaybreaks
\begin{align}\label{e:3.10}       
&R(P,Q,Q,P)\nonumber\\
&={\olr}(P,Q,Q,P)  
{-}\bigg(\!{g^\prime(r)\over g(r)}(A,A){-}{Lr\over 1{-}Lr^2}(X,X)\!
 \bigg)
 \bigg(\!{g^\prime(r)\over g(r)}(B,B){-}{Lr\over 1{-}Lr^2}(Y,Y)\!
 \bigg)\nonumber\\
&\qquad+\bigg({g^\prime(r)\over g(r)}(A,B)-{Lr\over 1-Lr^2}(X,Y)
        \bigg)^{\!2}\nonumber\\
&={\olr}(P,Q,Q,P)
-\bigg({g^\prime(r)\over g(r)}
 \bigg)^{\!2}\big((A,A)(B,B)-(A,B)^2
         \big)\nonumber\\
&\qquad-\bigg({Lr\over 1-Lr^2}
        \bigg)^{\!2}\big((X,X)(Y,Y)-(X,Y)^2
                \big)\nonumber\\
&\qquad+{g^\prime(r)\over g(r)}{Lr\over 1-Lr^2}
\big((A,A)(Y,Y)+(B,B)(X,X)-2(A,B)(X,Y)
\big).
\end{align}

To obtain the formula for the curvature at the points of~$M^\circ$, we
have to insert~\eqref{e:3.8}--\eqref{e:3.10} into~\eqref{e:3.7}.

\begin{lemma}\label{l:3.1}                
Let $p\in N\subset M$. Consider two vectors~$E$ and $F$ at the point~$p$ such 
that~$E=A+X$ and $F=B+Y$, where the vectors~$A$ and~$B$ are tangent to the fiber 
and the vector fields~$X$ and~$Y$ are basic. In this case the following formula
holds at the point~$p$:

\begin{align}\label{e:3.11}               
R(E,F,F,E)
&=R(X,Y,Y,X)-g^{\prime\prime\prime}(0)\big((A,A)(B,B)-(A,B)^2
                                      \big)\nonumber\\
&\qquad+L\big((A,A)(Y,Y)+(B,B)(X,X)-2(A,B)(X,Y)
         \big)\nonumber\\
&\qquad+3\big(B,\nabla_A[X,Y]
         \big).
\end{align}

\end{lemma}

\begin{proof}
Since the expressions on the right- and on the left-hand sides depend
on the values of the fields~$A$ and $B$ only at the point~$p$, we can
assume without loss of generality that the fields~$A$ and $B$ are
constant in the Cartesian coordinate system of the fiber. If $A=B=0$,
then~\eqref{e:3.11} is obvious. Therefore, henceforth we assume
that~$A\neq 0$. Let $a=\sqrt{(A,A)_p}$ and $b=(A,B)_p/a$. Consider the
ray emanating from~$p$ in the direction of the vector~$A$ and analyze
the limits of some expressions on that ray as~$r\to 0$. It is clear
that $A=a\partial/\partial r$ on the ray. Next, consider the field
$\OB=B-b\partial/\partial r$ on the ray. In~that case~$\OB$ is tangent
to the radial fibers. Applying the formulas for curvature at the
regular points of the ray obtained before, we get $(C=a\OB$ and
$Z=aY-bX)$

\begin{align*}
&R\bigg(a\pop r+X,b\pop r+\OB+Y,b\pop r+
                          \OB+Y,a\pop r+X
  \bigg)\\
&\qquad=R(X,\OB{+}Y,
      \OB{+}Y,X){-}{g^{\prime\prime}(r)\over g(r)}(C,C){+}
      \bigg({L\over 1{-}Lr^2}{+}{(Lr)^2\over(1{-}Lr^2)^2}
      \bigg)(Z,Z)\\
&\qquad\qquad+3\bigg({Lr\over 1-Lr^2}+{g^\prime(r)\over g(r)}
               \bigg)\big(C,[X,Y]
                     \big).
\end{align*}

But the value of $R(E,F,F,E)$ at the point~$p$ is the limit of this
expression as $r\to 0$ and equals

\begin{equation}\label{e:3.12}                  
R(X,{\OB}+Y,{\OB}+Y,X)-a^2g^{\prime\prime\prime}(0)(\OB,\OB)+L(Z,Z)+3A
\big(\OB,[X,Y]
\big).
\end{equation}

Note that $\big(\partial/\partial r,[X,Y]\big)=0$. Therefore,
$$
A\big(\OB,[X,Y]
 \big)=A\big(B,[X,Y]
        \big).
$$
Using the symmetries of the curvature tensor, from~\eqref{e:3.12} we
obtain the following equality at the point~$p$:

\begin{equation}\label{e:3.13}                  
R(X,\OB+Y,\OB+Y,X)=R(X,Y,Y,X)+L(\OB,\OB)(X,X).
\end{equation}

It remains to put~\eqref{e:3.13} into~\eqref{e:3.12} and use the equality
$(\nabla_AB)_p=\nobreak0$.
\end{proof}

Consider vectors~$A$, $B$, $X$, and~$Y$ in some vector space with an inner 
product~$\langle\boldsymbol{\cdot}\,,\boldsymbol{\cdot}\rangle$. Define

\begin{gather}
V(A,B)=\langle A,A
       \rangle
       \langle B,B
       \rangle-
       \langle A,B
       \rangle^2,\label{e:3.14}\\    
W(A,B;X,Y)=\langle A,A
           \rangle
           \langle Y,Y
           \rangle+
           \langle B,B
           \rangle\langle X,X
           \rangle-2\langle A,B
           \rangle\langle X,Y
           \rangle.\label{e:3.15}        
\end{gather}

\begin{lemma}\label{l:3.2}               
The following statements are true:
\begin{enumerate}
\item $W(A,B;X,Y)\geq 0$.
\item If $\{A,B\}\perp\{X,Y\}$ and the vectors $A+X$ and $B+Y$ form an 
orthonormal system, then

\begin{equation}\label{e:3.16}           
W(A,B;X,Y)\geq{1\over 2}
\big(\langle A,A\rangle+\langle B,B\rangle
\big)\big(\langle X,X
          \rangle+\langle Y,Y
                  \rangle
     \big).
\end{equation}

\item If $\{A,B\}\perp\{X,Y\}$, then

\begin{equation}\label{e:3.17}           
V(A+X,B+Y)=V(A,B)+V(X,Y)+W(A,B;X,Y).
\end{equation}

\end{enumerate}
\end{lemma}

\begin{proof}
1. It follows from the Cauchy--Bunyakovski\u\i{} inequality that

\begin{align*}
W(A,B;X,Y)
  &\geq\langle A,A\rangle
       \langle Y,Y\rangle+
       \langle B,B\rangle
       \langle X,X\rangle\\
  &\qquad-2\sqrt{\langle A,A\rangle
                 \langle B,B\rangle
                 \langle X,X\rangle
                 \langle Y,Y\rangle}\\
  &=\big(\sqrt{\langle A,A\rangle
               \langle Y,Y\rangle}-
         \sqrt{\langle B,B\rangle
               \langle X,X\rangle}\,
    \big)^2.
\end{align*}

2.
Consider the matrix
$$
G_{AB}=\begin{pmatrix}
\langle A,A\rangle &\langle A,B\rangle\\
\langle A,B\rangle &\langle B,B\rangle\\
\end{pmatrix}.
$$

Let $\lambda$ and $\mu$ be its eigenvalues; it is obvious that
$\lambda,\mu\in\mathbb R$. Furthermore,

\begin{align*}
W(A,B;X,Y)
&=\langle A,A\rangle
\big(1-\langle B,B\rangle
\big)+\langle B,B\rangle
\big(1-\langle A,A\rangle
\big)+2\langle A,B\rangle^2\\
&=\langle A,A\rangle+
  \langle B,B\rangle-2V(A,B)\\
&=\lambda+\mu-2\lambda\mu,
\end{align*}
\begin{align*}
\big(\langle A,A\rangle+
     \langle B,B\rangle
\big)\big(\langle X,X\rangle+
          \langle Y,Y\rangle
     \big)\\
&=(\lambda+\mu)(2-\lambda-\mu).
\end{align*}

Thus, the required inequality is equivalent to the following:
$$
\lambda+\mu-2\lambda\mu
\geq
\lambda+\mu-{1\over 2}(\lambda+\mu)^2.
$$
The last inequality is obvious.

3. This statement is verified by straightforward computations.
\end{proof}

\begin{lemma}\label{l:3.3}         
The following inequality holds for some constant $M_1$ depending only
on the metric $(\boldsymbol{\cdot}\,,\boldsymbol{\cdot})_0$:

\begin{equation}\label{e:3.18}      
\big\vert [X,Y]^\prime
\big\vert\leq M_1g(r){|X|\cdot |Y|\over 1-Lr^2}.
\end{equation}

Here the fields~$X$ and~$Y$ are basic.
\end{lemma}

\begin{proof}
Fix the fields on~$N$ that produce~$X$ and~$Y$. Then the expression on
the right-hand side is independent of~$L$ \big(since at a
point~$(p,v)$ it is equal to \hbox{$M_1g(r)|X_p|\cdot|Y_p|$}\big).
Moreover, the fields~$X$ and~$Y$ do not depend on~$L$ and~$g$; since
the orthogonal projection to the fibers of each fixed vector field is
also independent of~$L$ and~$g$, the field~$[X,Y]^\prime$ depends only
on the original metric~$(\boldsymbol{\cdot}\,,\boldsymbol{\cdot})_0$.
Therefore, we may assume that~$L=0$ and~$g(r)=r$. It remains to choose
the constant~$M_1$ such that the inequality~\eqref{e:3.18} holds in
that case. This is possible since the value~$[X,Y]^\prime$ at the
point~$(p,v)$ depends only on the values~$X_p$ and~$Y_p$ bilinearly
(and does not depend on the values of~$X$ and~$Y$ at other points);
moreover, at the points of~$N$ we have $[X,Y]^\prime=0$.
\end{proof}

\begin{lemma}\label{l:3.4}          
Let the metric~$g_N$ have positive sectional curvature and let
$g^{\prime\prime\prime}(0)<0$. Then there exists a constant~$L_1$ such
that the metric~$(\boldsymbol{\cdot}\,,\boldsymbol{\cdot})$ has
positive sectional curvature at the points of~$N$ provided that~$L\geq
L_1$.
\end{lemma}

\begin{proof}
Fix a point~$p\in N$. We may assume that the fields~$A$ and~$B$ are
constant in the Cartesian coordinate system of the fibers. Then
$$
\big(B,\nabla_A[X,Y]
\big)_p=A\big(B,[X,Y]
         \big)_p.
$$
Consider the ray emanating from~$p$ in the direction of the vector~$A$ and use the previous
lemma. We have

\begin{align*}
A\big(B,[X,Y]
 \big)_p
 &=|A|\cdot\pop r\biggr\vert_{r=0}\big(B,[X,Y]
                                  \big)\\
 &=|A|\lim_{r\to 0}{\big(B,[X,Y]
                    \big)\over r}\\
 &\leq |A|\cdot |B|\cdot M_1\cdot |X|\cdot |Y|\lim_{r\to 0}{g(r)\over r}\\
 &=M_1|A|\cdot|B|\cdot|X|\cdot|Y|.
\end{align*}

Therefore, we have the following inequality at the point~$p$:
$$
\big\vert\big(B,\nabla_A[X,Y]
         \big)
\big\vert\leq{M_1\over 4}
\big(|A|^2+|B|^2
\big)\big(|X|^2+|Y|^2
     \big).
$$
Choose $\delta>0$ such that $-g^{\prime\prime\prime}(0)\geq\delta$ and
$R(X,Y,Y,X)\geq\delta V(X,Y)$. Using Lemmas~\ref{l:3.1}
and~\ref{l:3.2}, we obtain the following inequality for an orthonormal
system~$E$ and~$F$:

\begin{align*}
R(E,F,F,E)
&\geq\delta\big(V(X,Y)+V(A,B)
           \big)+L\cdot W(A,B;X,Y)\\
&\qquad-2M_1\cdot W(A,B;X,Y)\\
&=\delta+(L-\delta-2M_1)\cdot W(A,B;X,Y).
\end{align*}

It remains to take $L\geq 2M_1+\delta$.
\end{proof}

\begin{lemma}\label{l:3.5}         
For various~$\varepsilon>0$, consider the
metrics~$(\boldsymbol{\cdot}\,,\boldsymbol{\cdot})$ constructed
using~\eqref{e:3.1} from the fixed
metric~$(\boldsymbol{\cdot}\,,\boldsymbol{\cdot})_0$ and the
function~$g_\varepsilon$ from Lemma~{\rm{\ref{l:1.2}}}. Suppose that
the metric~$g_N$ has positive sectional curvature. Then there exist
constants~$L_0$ and~$\rho_0>0$ independent of~$\varepsilon$ and such
that for~$L\geq L_0$ the
metric~$(\boldsymbol{\cdot}\,,\boldsymbol{\cdot})$ has positive
sectional curvature on~$\{r<\rho_0/L\}$.
\end{lemma}

\begin{proof}
The positivity of the sectional curvature at the points of~$N$ follows
from the previous lemma, so we have only to establish it at the points
of~$M^\circ$. For some constant~$\lambda\in\mathbb R$ consider the
metric~$(\boldsymbol{\cdot}\,,\boldsymbol{\cdot})_3$ constructed
using~\eqref{e:3.1} from the fixed
metric~$(\boldsymbol{\cdot}\,,\boldsymbol{\cdot})_0+\lambda\cdot\pi^*g_N$,
the function~$g=g_0$ from section~1 and a constant~$L_3$ so big that
this metric would have positive curvature at the points of~$N$
when~$\lambda=0$. In this case the sectional curvature of the
metric~$(\boldsymbol{\cdot}\,,\boldsymbol{\cdot})_3$ in some
neighborhood of~$N$ is bounded from below by some number~$\delta_3>0$
when~$\lambda$ lies in some neighborhood of zero.

Consider a radial fiber $\{r=t\}$. If~$t$ is sufficiently small, we can find
the value~$t_0$ such that $g_0(t_0)=g_\varepsilon(t)$.
The metric~$(\boldsymbol{\cdot}\,,\boldsymbol{\cdot})$ on the radial fiber~$\{r=t\}$
coincides with the metric~$(\boldsymbol{\cdot}\,,\boldsymbol{\cdot})_3$ on the radial fiber~$\{r=t_0\}$
for~$\lambda=(1-Lt^2)/(1-L_3t_0^2)-1$, the latter being close to~0 for sufficiently small~$t$
and~$Lt^2$. Then the curvature tensor~$\olr$ of the 
restriction of the metric~$(\boldsymbol{\cdot}\,,\boldsymbol{\cdot})$ to the radial fiber~$\{r=t\}$
coincides with the curvature tensor~${\olr}_3$ of the restriction of the
metric~$(\boldsymbol{\cdot}\,,\boldsymbol{\cdot})_3$ to the radial
fiber~$\{r=t_0\}$; thus, from~\eqref{e:3.10} we conclude that

\begin{align*}
R(P,Q,Q,P)
&\geq\delta_3 V(P,Q)+{g_0^\prime(t_0)^2-g_\varepsilon^\prime(t)^2
                      \over g_\varepsilon(t)^2}V(A,B)\\
&\qquad+\Bigg(\bigg({L_3t_0\over 1-L_3t_0^2}
              \bigg)^{\!2}-
              \bigg({Lt\over 1-Lt^2}
              \bigg)^{\!2}
        \Bigg)V(X,Y)\\
&\qquad+\bigg({g_\varepsilon^\prime(t)
               \over g_\varepsilon(t)}{Lt\over 1-Lt^2}-
               {g_0^\prime(t_0)\over g_0(t_0)}
               {L_3t_0\over 1-L_3t_0^2}
        \bigg)W(A,B;X,Y).
\end{align*}

Using Lemma~\ref{l:1.3}, we obtain the following inequality for
large~$L$ and small~$t$ and~$Lt$:
$$
R(P,Q,Q,P)\geq{\delta_3\over 2}V(P,Q)+{L\over 7}W(A,B;X,Y).
$$
It follows from here and Lemma~\ref{l:3.3} that

\begin{align*}
R(E,F,F,E)
&\geq {\delta_3\over 2} V(P,Q)+{L\over 7}W(A,B;X,Y)\\
\noalign{\vskip-3mm}
&\qquad+L(Z,Z)+(C,C)-4M_1|C|\cdot|X|\cdot|Y|.
\end{align*}

But
$$
W(A,B;X,Y)+(Z,Z)=W(\hA,\hB;X,Y);
$$
therefore,

\begin{align}\label{e:3.19}               
R(E,F,F,E)
&\geq {\delta_3\over 2}V(P,Q)+{L\over 7}W(\hA,\hB;X,Y)\nonumber\\
&\qquad+(C,C)-4M_1|C|\cdot|X|\cdot|Y|.
\end{align}

Assume that $W(\hA,\hB;X,Y)\leq\delta$ for some $\delta>0$. We have

\begin{align*}
&V(P,Q)+W(\hA,\hB;X,Y)+(C,C)\\
&\qquad\geq V(A,B)+(C,C)+V(X,Y)+W(\hA,\hB;X,Y)=1.
\end{align*}

On the other hand, it follows from Lemma~\ref{l:3.2}~(item~2) that the
value~$|C|\cdot|X|\cdot|Y|$ is infinitesimal as~$\delta\to 0$, so
for~$\delta$ sufficiently small we have~$R(E,F,F,E)>\nobreak0$. If, on
the contrary, $W(\hA,\hB;X,Y)\geq\delta$, then, recalling that all the
values participating in the inequality~\eqref{e:3.19} are bounded, we
can choose~$L$ such that ~$R(E,F,F,E)>0$ in this case as well. It
remains to choose~$\rho_0$ sufficiently small.
\end{proof}

\section{Gluing Metrics in a Neighborhood of a Submanifold}  

The proofs of the lemmas in this section base on some ideas taken
from~\cite[statement~2.1]{Ga}.

\begin{lemma}\label{l:4.1}       
Let $\varepsilon>0$. Then there exists a function
$\varphi:(0,+\infty)\to\mathbb R$ with the following properties:
\begin{enumerate}
\item $\varphi\in C^\infty$.
\item $0\leq\varphi(x)\leq 1$ for all $x$.
\item $\varphi(x)=1$ for $x\leq\delta_1$ and $\varphi(x)=0$ for
$x\geq\delta_2$, where $\delta_1$ and $\delta_2$ are some constants such that
$0<\delta_1<\delta_2\leq\varepsilon$.
\item $\big\vert x\varphi^\prime(x)\big\vert\leq\varepsilon$,
$\big\vert x^2\varphi^{\prime\prime}(x)\big\vert\leq\varepsilon$ for all~$x$.
\end{enumerate}
\end{lemma}

\begin{proof}
Consider a function~$f\in C^\infty(\mathbb R)$ such that $0\leq f\leq 1$,
$f(x)=1$ for~$x\leq 1$, and~$f(x)=0$ for~$x\geq 2$. Choose a constant~$N$
such that the following inequalities hold for all~$x$:
$$
\big\vert f^\prime(x)
\big\vert\leq N,
\quad
\big\vert f^{\prime\prime}(x)
\big\vert\leq N.
$$
Take $\delta>0$ and $\lambda\in (0,1)$ and put
$$
\varphi(x)=f\Big({x^\lambda\over\delta}
            \Big).
$$
It is clear that, for the fixed $\lambda$, one can choose~$\delta$
such that the function~$\varphi$ would satisfy the third condition.
Moreover, it is obvious that the function~$\varphi$ satisfies the
first two conditions. We will now show that one can choose~$\lambda$
such that for each~$\delta$ the function~$\varphi$ would satisfy the
fourth condition. We have

\begin{gather*}
\varphi^\prime(x)=f^\prime
\bigg({x^\lambda\over\delta}
\bigg)\delta^{-1}\lambda x^{\lambda-1},\\
\varphi^{\prime\prime}(x)=f^{\prime\prime}
\bigg({x^\lambda\over\delta}
\bigg)\delta^{-2}\lambda^2 x^{2\lambda-2}+f^\prime
\bigg({x^\lambda\over\delta}
\bigg)\delta^{-1}\lambda(\lambda-1)x^{\lambda-2}.
\end{gather*}

Therefore, it will suffice to satisfy the following inequalities:

\begin{gather*}
N\delta^{-1}\lambda x^{\lambda}\leq\varepsilon,\\
N(\delta^{-2}\lambda^2x^{2\lambda}+\delta^{-1}\lambda(1-\lambda)x^\lambda)
\leq\varepsilon.
\end{gather*}

But we have~$\varphi(x)=0$ for~$x^\lambda>2\delta$, so we need to
check the inequalities only in the case $x^\lambda\leq 2\delta$; thus,
it is enough to obtain the inequalities

\begin{gather*}
N\delta^{-1}\lambda 2\delta\leq\varepsilon,\\
N(\delta^{-2}\lambda^2(2\delta)^2+\delta^{-1}\lambda(1-\lambda)2\delta)
\leq\varepsilon\\
\end{gather*}

or, equivalently,

\begin{gather*}
2N\lambda\leq\varepsilon,
\quad
2N\lambda(1+\lambda)\leq\varepsilon.
\end{gather*}

Since the left-hand sides of the inequalities converge to zero
as~$\lambda\to 0$, we can choose~$\lambda$ which satisfies the
required conditions.
\end{proof}

\begin{lemma}\label{l:4.2}               
Define the function $\psi:\mathbb R^n\to\mathbb R$ as follows:
$\psi(x)=\varphi\big(|x|\big)$. Then
$|x|\cdot\big\vert\nabla\psi(x)\big\vert\leq\varepsilon$ and
$|x|^2\cdot\big\vert\nabla^2\psi(x)\big\vert\leq 2\varepsilon$ for all~$x$.
\end{lemma}

\begin{proof}
We have

\begin{gather*}
{\partial\psi\over\partial x_i}=
\varphi^\prime\big(|x|
              \big){x_i\over|x|},
\quad
{\partial^2\psi\over\partial x_i\partial x_j}=
\varphi^{\prime\prime}\big(|x|
                      \big){x_ix_j\over|x|^2}+\varphi^\prime
                      \big(|x|
                      \big){\delta_{ij}|x|^2-x_ix_j\over|x|^3};
\end{gather*}

therefore,

\begin{gather*}
|x|\cdot\left\vert\partial\psi\over\partial x_i
        \right\vert\leq\varepsilon,
\quad
|x|^2\cdot
\left\vert\partial^2\psi\over\partial x_i\partial x_j
\right\vert\leq 2\varepsilon.
\end{gather*}

\end{proof}

Suppose that $M$ is a Riemannian manifold with the metric~$g_0$, $N$
is its compact submanifold, $m=\dim M$, and $n=\dim N$. Define the
map~\hbox{$r:M\to\mathbb R$} using the following formula:

\begin{equation}\label{e:4.1}      
r(\exp_pv)=|v|.
\end{equation}

Here $p\in N$ and $v\in T_pN^\perp$.

\begin{lemma}\label{l:4.3}           
Assume that we have a metric~$g_1$ in a neighborhood of~$N$ such that
\begin{enumerate}
\item $g_0$ and $g_1$ have positive sectional curvature.
\item The $1$-jets of two metrics coincide at the points of~$N$.
\end{enumerate}

Then, for each $\rho>0$, there exists a smooth Riemannian metric~$g$
on~$M$ having positive sectional curvature and such that~$g=g_0$
for~$r>\rho$ and $g=g_1$ for~$r$ sufficiently small.
\end{lemma}

\begin{proof}
In a neighborhood of~$N$, define the family of metrics

\begin{equation}\label{e:4.2}          
\overline{g}(q,s)=(1-s)g_0(q)+sg_1(q),\ \ s\in[0,1].
\end{equation}

Choose the function~$\varphi$ from Lemma~\ref{l:4.1} for
some~$\varepsilon\leq\rho$. Put

\begin{equation}\label{e:4.3}           
g(q)=\overline{g}
\Big(q,\varphi
     \big(r(q)
     \big)
\Big).
\end{equation}

It follows from the properties of the function~$\varphi$ that $g$~is a
Riemannian metric which coincides with~$g_0$ for~$r>\rho$ and
with~$g_1$ for a sufficiently small~$r$. We will show that, for each
point~$p\in N$, there are an open neighborhood~$U(p)$ and a
number~$\varepsilon(p)>0$ such that the metric~$g$ has positive
curvature in~$U(p)$ for $\varepsilon\leq\varepsilon(p)$. Having
obtained that, we use the compactness of~$N$ to conclude that for
sufficiently small~$\varepsilon$, the metric~$g$ has positive
sectional curvature in some neighborhood~$U$ of the submanifold~$N$
independent of~$\varepsilon$. We are left with taking~$\varepsilon$ so
small that~$g$ would have positive curvature in~$U$ and the
set~$\{r\leq\varepsilon\}$ would lie inside~$U$. Then~$g$ coincides
with~$g_0$ outside~$U$ and therefore has positive curvature
outside~$U$. Thus $g$ is the sought metric.

Fix $p\in N$ and consider a coordinate system $(x_1,\dots,x_m)$ in a
neighborhood of~$p$ such that~$x_i(p)=0$ and
$r^2=x_1^2+\dots+x_{m-n}^2$. Such system can be constructed by any
trivialization of the normal bundle over~$N$. It is clear that in this
case the equations from Lemma~\ref{l:4.2} hold for the function
$\psi(x)=\varphi\big(r(x)\big)$.

We use the following formulas for the sectional curvature:

\begin{gather*}
R(X,Y,Y,X)=R_{ijkl}X^i Y^j Y^k X^l,
\quad
R_{ijkl}=R_{jkl}^\alpha g_{\alpha i},\\
R_{ijk}^l={\partial\Gamma_{ik}^l\over\partial x^j}-
{\partial\Gamma_{ij}^l\over\partial x^k}+
\Gamma_{j\alpha}^l\Gamma_{ik}^\alpha-
\Gamma_{k\alpha}^l\Gamma_{ij}^\alpha,\\
\Gamma_{ij}^k={1\over 2}g^{kl}
\left({\partial g_{il}\over\partial x^j}+
      {\partial g_{jl}\over\partial x^i}-
      {\partial g_{ij}\over\partial x^l}
\right).
\end{gather*}

Denote by $\Gamma$ and $R$ the Christoffel symbols and the components
of the curvature tensor of the metric~$g$. Also, introduce the
notations $\overline\Gamma{}_{ij}^{\,k}(s)$ and $\olr_{ijk}^{\,l}(s)$
for the Christoffel symbols and the components of the curvature tensor
of the metric~$\overline g$ for some fixed value of the parameter~$s$.
Denote by dot the derivative with respect to~$s$.
Note that
$$
\dot{\overline g}_{ij}=g_{1ij}-g_{0ij},
\quad
\ddot{\overline g}_{ij}=0.
$$
Since the metrics~$g_0$ and~$g_1$ together with their first
derivatives coincide at the points of~$N$, we have the following
relations at~$r=0$:
$$
\dot{\overline{g}}_{ij}=0,
\quad
\dot{\overline{g}}{}^{\,ij}=0,
\quad
{\partial\dot{\overline{g}}_{ij}\over\partial x^k}=0,
\quad
{\partial\dot{\overline{g}}{}^{\,ij}\over\partial x^k}=0,
\quad
\dot{\overline\Gamma}{}_{ij}^{\,k}=0,
\quad
\ddot{{\olr}}_{ijkl}=0.
$$
Using the Taylor formula, we conclude that the following inequalities
hold for some constant~$N$ and sufficiently small~$r$:

\begin{alignat*}3
&\big\vert g_{ij}
 \big\vert\leq N,
  &&\quad\big\vert g^{ij}
         \big\vert\leq N,
    &&\quad\bigg\vert{\partial g^{ij}\over\partial x^k}
           \bigg\vert\leq N,\\
&\big\vert\dot{\overline{g}}_{ij}
 \big\vert\leq Nr^2,
  &&\quad\big\vert\dot{\overline{g}}{}^{\,ij}
         \big\vert\leq Nr^2,
    &&\quad\bigg\vert{\partial\dot{\overline g}_{ij}\over\partial x^k}
           \bigg\vert\leq Nr,\\
&\big\vert\Gamma_{ij}^k
 \big\vert\leq N,
  &&\quad\big\vert\overline{\Gamma}{}_{ij}^k(s)
         \big\vert\leq N,
    &&\quad\big\vert\dot{\overline \Gamma}{}_{ij}^k
           \big\vert\leq Nr.
\end{alignat*}

We have
$$
{\partial g_{ij}\over\partial x^k}=
{\partial\overline{g}_{ij}\over\partial x^k}
\big(\psi(x)
\big)+\dot{\overline g}_{ij}
{\partial \psi\over\partial x^k},
\quad
\Gamma_{ij}^k=\overline\Gamma{}_{ij}^{\,k}
\big(\psi(x)
\big)+\Delta\Gamma_{ij}^k,
$$
where
$$
\Delta\Gamma_{ij}^k={1\over 2}g^{kl}
\left(\dot{\overline g}_{il}{\partial\psi\over\partial x^j}+
      \dot{\overline g}_{jl}{\partial\psi\over\partial x^i}-
      \dot{\overline g}_{ij}{\partial\psi\over\partial x^l}
\right);
$$
moreover,

\begin{gather*}
{\partial\Gamma_{ij}^k\over\partial x^l}=
{\partial\overline\Gamma{}_{ij}^{\,k}\over\partial x^l}\big(\psi(x)
                                                       \big)+
\dot{\overline\Gamma}{}_{ij}^{\,k}\big(\psi(x)
                                   \big)
{\partial\psi\over\partial x^l}+
{\partial\Delta\Gamma_{ij}^k\over\partial x^l},\\
{\partial\over\partial x^\alpha}
\left(g^{kl}\dot{\overline g}_{il}{\partial\psi\over\partial x^j}
\right)=
{\partial g^{kl}\over\partial x^\alpha}\dot{\overline g}_{il}
{\partial\psi\over\partial x^j}+g^{kl}
\left({\partial\dot{\overline g}_{il}\over\partial x^\alpha}
      {\partial\psi\over\partial x^j}+\dot{\overline g}_{il}
      {\partial^2\psi\over\partial x^jx^\alpha}
\right),\\
R_{ijkl}={\olr}{}_{ijkl}\big(\psi(x)
                        \big)+g_{\alpha i}\Delta R_{jkl}^\alpha,
\end{gather*}
where
\begin{align*}
\Delta R_{ijk}^l
&=\dot{\overline\Gamma}{}_{ik}^{\,l}\big(\psi(x)
                                    \big)
  {\partial\psi\over\partial x^j}+
  {\partial\Delta\Gamma_{ik}^{l}\over\partial x^j}-
  \dot{\overline\Gamma}{}_{ij}^{\,l}\big(\psi(x)
                                    \big)
  {\partial\psi\over\partial x^k}-
  {\partial\Delta\Gamma_{ij}^{l}\over\partial x^k}\\
&\qquad+\overline{\Gamma}{}_{j\alpha}^{\,l}\big(\psi(x)
                                           \big)\Delta\Gamma_{ik}^\alpha+
  \Gamma{}_{ik}^\alpha\Delta\Gamma_{j\alpha}^{\,l}-
  \overline{\Gamma}{}_{k\alpha}^{\,l}\big(\psi(x)
                                     \big)\Delta\Gamma_{ij}^\alpha-
  \Gamma{}_{ij}^\alpha\Delta\Gamma_{k\alpha}^l.
\end{align*}

Estimate the values $\Delta\Gamma_{ij}^k$,
$\partial\Delta\Gamma_{ij}^k\over\partial x^l$, and $g_{\alpha
i}\Delta R_{jkl}^\alpha$ as follows:

\begin{align*}
|\Delta\Gamma_{ij}^k|
&\leq 2N^2mr^2\big\vert\nabla\psi(x)
              \big\vert\leq 2N^2m\varepsilon r,\\
\bigg\vert{\partial\Delta\Gamma_{ij}^k\over\partial x^l}
\bigg\vert
&\leq 2N^2mr^2\big\vert\nabla\psi(x)
              \big\vert+2N^2mr\cdot
              \big\vert\nabla\psi(x)
              \big\vert+2N^2mr^2\big\vert\nabla^2\psi(x)
              \big\vert\\
&\leq 2N^2m\varepsilon(3+r),\\
|g_{\alpha i}\Delta R_{jkl}^\alpha|
&\leq 2N^2mr\cdot\big\vert\nabla\psi(x)
                 \big\vert+4N^3m^2\varepsilon(3+r)+8N^4m^3\varepsilon r\\
&\leq 2N^2m\varepsilon\big(1+2Nm(3+r)+4N^2m^2r
                      \big).
\end{align*}

As we see, with a proper choice of $\varepsilon$ one can make the
value $R_{ijkl}-{\olr}_{ijkl}\big(\psi(x)\big)$ arbitrarily small. It
remains to prove that we have the following inequality in a
neighborhood of~$p$ for some constant~$C>0$, all~$s\in [0,1]$ and
every orthonormal system~$X$ and~$Y$:
$$
{\olr}(s)(X,Y,Y,X)>C.
$$

Since $\ddot{\olr}_{ijkl}=0$ at the points of~$N$, the following
relation holds at these points:

\begin{equation}\label{e:4.4}                   
{\olr}_{ijkl}(s)=(1-s)R_{0ijkl}+sR_{1ijkl},
\end{equation}

That is, the tensor ${\olr}$ is a convex combination of the curvature
tensors of the metrics~$g_0$ and~$g_1$. Since the metrics~$g_0$
and~$g_1$ have positive curvature, we can obtain the inequality
$$
{\olr}(s)(p)(X,Y,Y,X)>2C
$$
for some constant $C>0$. But the expression~${\olr}(s)(X,Y,Y,X)$ is a
continuous function of the variables~$x$, $s$, $X$, and~$Y$, the last
three taking values in some compact. Then this expression is uniformly
continuous for small~$x$; therefore, if~$x$ is sufficiently close
to~$p$, the following inequality holds:
$$
\big\vert{\olr}(s)(x)(X,Y,Y,X)-{\olr}(s)(p)(X,Y,Y,X)
\big\vert\leq C.
$$
From here we obtain the sought inequality for the curvature
tensor~${\olr}$ in a small neighborhood of the point~$p$. The proof is
complete.
\end{proof}

\section{Taking Quotients by a Nonfree Isometric Action}

In this section we prove the theorem stated in the introduction.

We only give the proof for the case~$G=S^3$ (the case $G=S^1$ is
analyzed similarly). Put $M^\circ=M\setminus N$. It follows from the
second condition of the theorem that the action is free on~$M^\circ$.
So the space~$M^\circ/G$ has the structure of a smooth manifold, and
the quotient metric~$(\boldsymbol{\cdot}\,,\boldsymbol{\cdot})_{0*}$
of the metric~$(\boldsymbol{\cdot}\,,\boldsymbol{\cdot})_0$ has
positive sectional curvature on it due to the first condition of the
theorem. On the other hand, some neighborhood of~$N$ in the
manifold~$M$ is diffeomorphic to some neighborhood of~$N$ in the
normal bundle over~$N$; the diffeomorphism can be constructed
using~\eqref{e:2.11}, namely $\Phi(p,v)=\exp_pv$, where $p\in N$ and
$v\in T_pN^\perp$. So henceforth we assume that~$M$ is a neighborhood
of~$N$ in the normal bundle over $N$.

Note that, since $N$ is the fixed point set of a group of isometries, it
is totally geodesic. Furthermore, since the orbit of~$G$ at each point
of~$N$ is a singleton, it follows from the first condition of the
theorem that the restriction of the
metric~$(\boldsymbol{\cdot}\,,\boldsymbol{\cdot})_0$ to~$N$ has
positive sectional curvature. As above, we denote this restriction
by~$g_N$.

Using~\eqref{e:3.1}, for some constant~$L$ construct a $g_0$-proper
metric~$(\boldsymbol{\cdot}\,,\boldsymbol{\cdot})_1$ in a neighborhood
of~$N$ such that its 1-jets coincide at the points of~$N$ with those
of the original metric~$(\boldsymbol{\cdot}\,,\boldsymbol{\cdot})_0$.
(Henceforth we use the functions~$g_0$ and~$g_\varepsilon$
from~\S\,1.) Since the
metric~$(\boldsymbol{\cdot}\,,\boldsymbol{\cdot})_1$ is unique, it is
$G$ invariant. Now, take~$L$ so large that Lemma~\ref{l:3.4}
guarantees the positivity of the sectional curvature of the
metric~$(\boldsymbol{\cdot}\,,\boldsymbol{\cdot})_1$ at the points
of~$N$. By narrowing the manifold~$M$ to a smaller neighborhood of~$N$
in the normal bundle if necessary, we can assume that the
metrics~$(\boldsymbol{\cdot}\,,\boldsymbol{\cdot})_0$
and~$(\boldsymbol{\cdot}\,,\boldsymbol{\cdot})_1$ have positive
curvature everywhere on~$M$. Denote the quotient metric of the
metric~$(\boldsymbol{\cdot}\,,\boldsymbol{\cdot})_1$
by~$(\boldsymbol{\cdot}\,,\boldsymbol{\cdot})_{1*}$.

Now we study the action of~$G$. For each point~$p\in N$ consider the
action of~$G$ on the space $V_p=T_pN^\perp$. It is a representation
of~$G$ on the fiber~$V_p$, which we denote by~$G_p$. If $g\in G$, then
denote by~$g_p$ the image of~$g$ under the representation~$G_p$.
Since~$G$ acts by isometries, we have the following equality for $g\in
G$:
$$
g\big(\Phi(p,v)
 \big)=\Phi(p,g_pv).
$$
Since~$M$ was assumed to be a part of the normal bundle, i.e., $\Phi$
was the identity map, the group~$G$ acts on each fiber~$V_p$ according
to the representation~$G_p$. Note that each of the
representations~$G_p$ has no nonzero nonfree points due to the second
condition of the theorem; therefore, it is equivalent to the
representation of~$Sp(1)$ in~$\mathbb H^2$ given by the formula
$g(q_1,q_2)=(g\cdot q_1,g\cdot q_2)$ (such representation will be
henceforth called {\it standard}).

Let $(U,\Theta)$ be some trivialization of~$M$. We may assume
that~$\Theta$ acts from $U\times\mathbb H^2$. From now on we consider
only such trivializations for which the action of~$G$ on~$M$
corresponds to the standard action on~$\mathbb H^2$, i.e., for $p\in
U$, $x\in\mathbb H^2$, and $g\in G$ holds
$$
g\Theta(p,x)=\Theta(p,g\cdot x).
$$
Such trivializations still exist in a neighborhood of every point
of~$N$. The map~$\Psi_{12}$ in~\eqref{e:2.4} now belongs to
$Sp(2)\subset SO(8)$; here the elements of~$Sp(2)$ are embedded as
right multiplications.
                                                                  
Next, we study the structure of the space~$\mathbb H^2/Sp(1)$. First
consider the space~$S^7/Sp(1)$, where $S^7\subset\mathbb H^2$ carries
the standard metric of radius one. Note that there is an isometric
right action of~$Sp(2)$ on~$S^7$ which commutes with the isometric
left action of~$Sp(1)$; therefore, $Sp(2)$ acts isometrically
on~$S^7/Sp(1)$. The kernel of the latter action is~$\mathbb Z_2=\{\pm
E\}$, so the isometry group of~$S^7/Sp(1)$ has dimension at
least~$\dim Sp(2)=10$. On the other hand, $\dim\big(S^7/Sp(1)\big)=4$.
As is well-known (see~\cite[Chapter~2, Theorem~3.1]{Ko}), in this case
we have~$S^7/Sp(1)\simeq S^4$ endowed with the standard metric
multiplied by some constant.

Consider the geodesic $c(\varphi)=(\cos\varphi,\sin\varphi)\in
S^7\subset\mathbb H^2$ for $\varphi\in [0,\pi]$. It is everywhere
orthogonal to the orbits of~$Sp(1)$, so its projection~$c_*$ to~$S^4$
is a geodesic of the same length~$\pi$. On the other hand,
$c_*(0)=c_*(\pi)$, and~$c_*$ has no other self-intersections. It
follows from the behavior of the geodesics on round sphere that the
metric on $S^4\simeq S^7/Sp(1)$ is the same as on the standard sphere
of radius~$1/2$.

Let $f:S^7\to S^4$ be the quotient map. Extend~$f$ to a map from~$\mathbb H^2$ 
to~$\mathbb R^5$ by the rule $f(tv)=tf(v)$,
where~$t\geq 0$ and $v\in S^7$. Note that the map~$f$ is everywhere continuous
and smooth everywhere except zero. The right action of~$Sp(2)$ on~$S^7/Sp(1)$
determines the group homomorphism~$\gamma:Sp(2)\to SO(5)$.
Let $\gamma_*:\mathfrak{sp}(2)\to\mathfrak{so}(5)$ be the corresponding
homomorphism of Lie algebras. Note that 
\begin{align*}
f(x\cdot A)   &=f(x)\cdot\gamma(A),\\
df_x(x\cdot B)&=f(x)\cdot\gamma_*(B)
\end{align*}
for each~$A\in Sp(2)$, $B\in\mathfrak{sp}(2)$, and $x\in\mathbb H^2$.

Let $(U,\Theta)$ be a trivialization. It naturally engenders
a homeomorphism $\Theta_*:U\times\mathbb R^5\to\pi^{-1}(U)/G$
such that
$$
\Theta_*\big(p,f(x)
        \big)=\tau\big(\Theta(p,x)
                  \big)
$$
for $p\in U$ and $x\in\mathbb H^2$.

Now, the smooth structure on~$M/G$ is constructed by taking as coordinate charts
the maps~$\Theta_*$ for various trivializations. Smoothness of~$\tau$
on~$M^\circ$ follows from smoothness of~$f$ everywhere except zero. We now show
smoothness of the transition maps. Let~$(U,\Theta_1)$ and
$(U,\Theta_2)$ be two trivializations and $\Psi_{12}$ the map that connects them
as in~\eqref{e:2.4}. Then
\vspace{-3mm}

\begin{align*}
\Theta_{2*}\big(p,f(x)
           \big)
&=\tau\big(\Theta_2(p,x)
      \big)\\
&=\tau\Big(\Theta_1\big(p,x\cdot\Psi_{12}(p)
                   \big)
      \Big)\\
&=\Theta_{1*}
      \Big(p,f\big(x\cdot\Psi_{12}(p)
              \big)
      \Big)\\
&=\Theta_{1*}\Big(p,f(x)\cdot\gamma\big(\Psi_{12}(p)
                                   \big)
             \Big)
\end{align*}
for~$p\in U$ and $x\in\mathbb H^2$.

Note that $M/G$ is a vector bundle of rank 5 over~$N/G=N$ with
trivializations~$(U,\Theta_*)$.

We proceed to constructing the sought positively curved metric
on~$M/G$. First consider a trivialization~$(U,\Theta)$. Let
$Q_p:T_pN\to\mathfrak{so}(8)$ be a map used in construction of the
metric~$(\boldsymbol{\cdot}\,,\boldsymbol{\cdot})_1$ in this
trivialization. Since the
metric~$(\boldsymbol{\cdot}\,,\boldsymbol{\cdot})_1$ and locally
horizontal fields are invariant under the action of~$G$, the
orthogonal projection of every locally horizontal field is invariant
under the action of~$G$. Thus, $Q_p(X)\in\nobreak\mathfrak{sp}(2)$ for
each~$X\in T_pN$. In this case we can define the map
$Q_{*p}:T_pN\to\mathfrak{so}(5)$ by the formula
$$
Q_{*p}(X)=\gamma_*\big(Q_p(X)
                  \big).
$$

Using~\eqref{e:3.1} and Lemma~\ref{l:2.3}, for each $\varepsilon>0$
construct the
metric~$(\boldsymbol{\cdot}\,,\boldsymbol{\cdot})_\varepsilon$
on~$\pi^{-1}(U)/G$ from the metric~$g_N$ on~$N$, the
function~$g_\varepsilon$ from Lemma~\ref{l:1.2}, the map~$Q_{*p}$, and
the constant~$L$. The constructed metric is the unique
$g_\varepsilon$-proper metric on the bundle~$\pi^{-1}(U)/G$ such that
it coincides with~$(\boldsymbol{\cdot}\,,\boldsymbol{\cdot})_{1*}$ for
$r\geq\varepsilon$ and the restriction of the map~$\pi$ to every
radial fiber is a Riemannian submersion onto~$N$ endowed with the
metric~$(1-Lr^2)g_N$.

Indeed, the equality of the metrics on the fibers follows from the
fact that $g_\varepsilon=g_0/2$ for $r\geq\varepsilon$ and on the
spaces orthogonal to the fibers, from the fact that $\pi$ is a
Riemannian submersion. It remains to prove that the orthogonal
projection of each locally horizontal field to the fiberw is the same
in both metrics. Let $X$ be such a field, $p\in N$, and $x\in\mathbb
H^2$. Then the orthogonal projection of the field~$X$ to the fibers in
the metric~$(\boldsymbol{\cdot}\,,\boldsymbol{\cdot})_\varepsilon$ is

\begin{align*}
X^{\prime\varepsilon}_{\tau\big(\Theta(p,x)
                           \big)}
&=d\Theta_{*\big(p,f(x)
            \big)}\big(0,f(x)\cdot Q_{*p}(X_p)
                  \big)\\
&=d\Theta_{*\big(p,f(x)
            \big)}\Big(0,df_x\big(x\cdot Q_p(X_p)
                             \big)
                  \Big)\\
&=d\tau_{\Theta(p,x)}\Big(d\Theta_{(p,x)}\big(0,x\cdot Q_p(X_p)
                                         \big)
                     \Big)\\
&=d\tau_{\Theta(p,x)}\big(X^{\prime1}_{\Theta(p,x)}
                     \big).
\end{align*}

Here $X^{\prime1}$ is the orthogonal projection to the fibers with
respect to the metric~$(\boldsymbol{\cdot}\,,\boldsymbol{\cdot})_1$ of
the locally horizontal field~$X^\circ$ on~$M$ produced by the same
field on~$N$ as~$X$. Note that~$X$ and~$X^\circ$ are $\tau$ related.
The expression to the right equals the orthogonal projection of the
field~$X$ to the fibers in the
metric~$(\boldsymbol{\cdot}\,,\boldsymbol{\cdot})_{1*}$.

The uniqueness of the
metric~$(\boldsymbol{\cdot}\,,\boldsymbol{\cdot})_\varepsilon$ makes
it possible to define it on the whole~$M/G$. Note that 1-jets of the
obtained metric at the points of~$N$ are independent of~$\varepsilon$.
Choose constants~$L_0$ and~$\rho_0$ in Lemma~\ref{l:3.5} such that the
metric~$(\boldsymbol{\cdot}\,,\boldsymbol{\cdot})_\varepsilon$ has
positive curvature for $L\geq L_0$, $r<\rho_0/L$, and
every~$\varepsilon$. Fix~$L\geq L_0$. It follows from
lemma~\ref{l:4.3} that there exists a positively curved
metric~$(\boldsymbol{\cdot}\,,\boldsymbol{\cdot})^\prime$ on~$M$ such
that $(\boldsymbol{\cdot}\,,\boldsymbol{\cdot})^\prime=
(\boldsymbol{\cdot}\,,\boldsymbol{\cdot})_0$ outside any neighborhood
of~$N$ given a priori,
and~$(\boldsymbol{\cdot}\,,\boldsymbol{\cdot})^\prime=
(\boldsymbol{\cdot}\,,\boldsymbol{\cdot})_1$ when~$r<2\varepsilon$ for
some~$\varepsilon>0$. We may assume that $2\varepsilon<\rho_0/L$. As
follows from the proof of Lemma~\ref{l:4.3} (according
to~\eqref{e:4.3} this metric is a linear combination of two $G$
invariant metrics with a coefficient that is constant on the orbits
of~$G$), the metric~$(\boldsymbol{\cdot}\,,\boldsymbol{\cdot})^\prime$
is invariant under the action of~$G$.

Let $(\boldsymbol{\cdot}\,,\boldsymbol{\cdot})^\prime_*$ be the
quotient metric of the metric
$(\boldsymbol{\cdot}\,,\boldsymbol{\cdot})^\prime$ on $M^\circ/G$.
For $\varepsilon<r<2\varepsilon$, we have
$(\boldsymbol{\cdot}\,,\boldsymbol{\cdot})^\prime_*=
(\boldsymbol{\cdot}\,,\boldsymbol{\cdot})_{1*}=
(\boldsymbol{\cdot}\,,\boldsymbol{\cdot})_\varepsilon$. Therefore,
there is a $C^2$ smooth metric
$(\boldsymbol{\cdot}\,,\boldsymbol{\cdot})_*$ on $M/G$ such that
$(\boldsymbol{\cdot}\,,\boldsymbol{\cdot})_*=
(\boldsymbol{\cdot}\,,\boldsymbol{\cdot})^\prime_*$ for
$r>\varepsilon$ and $(\boldsymbol{\cdot}\,,\boldsymbol{\cdot})_*=
(\boldsymbol{\cdot}\,,\boldsymbol{\cdot})_\varepsilon$ for
$r<2\varepsilon$. This metric coincides with the quotient metric of
the metric $(\boldsymbol{\cdot}\,,\boldsymbol{\cdot})_0$ outside a
neighborhood of~$N$. It has positive curvature for~$r>\varepsilon$ by
Lemma~\ref{l:4.3}, and for~$r<2\varepsilon$ by Lemma~\ref{l:3.5}.
Therefore $(\boldsymbol{\cdot}\,,\boldsymbol{\cdot})_*$ is the sought
metric; the theorem is proven.
\qed\medskip

\noindent{\sc Remark}

Unfortunately, the author has not succeeded in using the method described above for constructing new
examples of positively curved spaces. However, we managed to obtain several 
already known examples; two of them are described below.

1.
Let $M=\mathbb{CP}^3$; the group $G=S^1$ acts as follows:
$$
g\cdot(z_0:z_1:z_2:z_3)=(gz_0:gz_1:z_2:z_3).
$$
Here $g\in S^1\subset\mathbb C$ and $z_i\in\mathbb C$. The fixed point space~$N$
has two connected components:

\begin{align*}
N_1 &=\big\{(z_0:z_1:0:0)\bigm\vert(z_0:z_1)\in\mathbb{CP}^1
      \big\},\\
N_2 &=\big\{(0:0:z_2:z_3)\bigm\vert(z_2:z_3)\in\mathbb{CP}^1
      \big\}.
\end{align*}

Each of these components is diffeomorphic to~$\mathbb{CP}^1\simeq S^2$.
Therefore~$N$ has codimension~4. The space~$M/G$ is homeomorphic to the sphere~$S^5$.

2. Let
$M=Sp(3)/Sp(1)^3$. Here~$Sp(1)^3$ is embedded in~$Sp(3)$ as matrixes
of the form~$\diag(q_1,q_2,q_3)$ for~$q_1,q_2,q_3\in Sp(1)$. We use
the positively curved metric on~$M$ constructed in~\cite{Wa}. The
group $G=S^3=Sp(1)$ acts as follows: $g\cdot
[A]=\big[\diag(g,1,1)\cdot A\big]$. Here $g\in Sp(1)\subset\mathbb H$,
and $[A]\in M$ is the equivalency class of the matrix~$A\in Sp(3)$.
The fixed point set~$N$ has three connected components:

\begin{gather*}
N_1=\Big\{\big[\diag(1,B)
          \big]\Bigm\vert B\in Sp(2)
    \Big\},\\
N_2=N_1\cdot\begin{pmatrix} 0&1&0\\
                            1&0&0\\
                            0&0&1
            \end{pmatrix},    
\qquad N_3=N_1\cdot\begin{pmatrix} 0&0&1\\
                            0&1&0\\
                            1&0&0
            \end{pmatrix}.
\end{gather*}

Each of the components is diffeomorphic to~$Sp(2)/Sp(1)^2\simeq S^4$.
Therefore $N$ has codimension~8. The space~$M/G$ is homeomorphic to
the sphere~$S^9$. \smallskip

The author is grateful to Ya.V.Baza\u\i kin for stating the problem
and helpful advice.

{\sc Novosibirsk State University, Novosibirsk, Russia}

{\it E-mail address:} {\tt dyatlov@gmail.com}

\end{document}